\documentclass[11pt]{amsart}
\usepackage{amsmath, amsthm, amssymb, latexsym,url}
\usepackage{hyperref}

\newtheorem{thm}{Theorem}[section]

\newtheorem{rem}[thm]{Remark}
\newtheorem{lem}[thm]{Lemma}

\newtheorem{prop}[thm]{Proposition}

\allowdisplaybreaks

\title{Towards a sharp converse of Wall's theorem on arithmetic progressions}

\author{Joseph Vandehey}

\begin{document}

\maketitle

\begin{abstract}
Wall's theorem on arithmetic progressions says that if $0.a_1a_2a_3\dots$ is normal, then for any $k,\ell\in \mathbb{N}$, $0.a_ka_{k+\ell}a_{k+2\ell}\dots$ is also normal. We examine a converse statement and show that if $0.a_{n_1}a_{n_2}a_{n_3}\dots$ is normal for periodic increasing sequences $n_1<n_2<n_3<\dots$ of asymptotic density arbitrarily close to $1$, then $0.a_1a_2a_3\dots$ is normal. We show this is close to sharp in the sense that there are numbers $0.a_1a_2a_3\dots$ that are not normal, but for which $0.a_{n_1}a_{n_2}a_{n_3}\dots$ is normal along a large collection of sequences whose density is bounded a little away from $1$.
\end{abstract}

\section{Introduction}

We will fix an integer base $b\ge 2$ throughout this paper.

Suppose $x\in [0,1)$ has (base-$b$) expansion $x=0.a_1a_2a_3\dots$. We say that $x$ is (base-$b$) normal if for every finite string $s=[d_1,d_2,\dots,d_k]$ with $d_i\in \{0,1,\dots,b-1\}$, we have that
\begin{equation}\label{eq:defn}
\lim_{n\to \infty} \frac{\#\{0\le i \le n-1: a_{i+j}=d_j,j=1,2,\dots,k\}}{n} = \frac{1}{b^k}.
\end{equation}
In other words, a number is normal if every string appears with the same limiting frequency as every other string of the same length.

In his thesis, Donald Dines Wall \cite{WallThesis} proved that selection along arithmetic progressions preseves normality. In other words, if $0.a_1a_2a_3\dots$ is normal, then for every $k,\ell \in \mathbb{N}$, $0.a_ka_{k+\ell}a_{k+2\ell}\dots$ is also normal. This we will refer to as Wall's theorem on arithmetic progressions.

At a recent conference on normal numbers in Vienna, Bill Mance described Wall's theorem on arithmetic progression as an ``if and only if" statement. That is, ``A number $0.a_1a_2a_3\dots$ is normal if and only if for every $k,\ell \in \mathbb{N}$, $0.a_ka_{k+\ell}a_{k+2\ell}\dots$ is normal." In the forward direction, this is just Wall's theorem as it is typically stated. In the reverse direction, this is trivial, since by letting $k=\ell=1$, then the number $0.a_ka_{k+\ell}a_{k+2\ell}\dots$ is just $0.a_1a_2a_3\dots$. Indeed, it can quickly be seen that for any $k\in \mathbb{N}$ and $\ell=1$, the normality of $0.a_ka_{k+\ell}a_{k+2\ell}\dots$ immediately gives the normality of $0.a_1a_2a_3\dots$.

However, it is reasonable to ask: if these trivial cases are removed, is Wall's theorem still an ``if and only if" statement?

We answer this in the negative.

\begin{thm}\label{thm:one}
There exists a real number $0.a_1a_2a_3\dots\in[0,1)$ that is \emph{not} normal such that for every $k\in \mathbb{N}$ and every $\ell\in \mathbb{N}$, $\ell\ge 2$, the number $0.a_ka_{k+\ell}a_{k+2\ell}\dots$ is normal.
\end{thm}

In particular, if the number $0.a_1a_2a_3\dots$ is normal, then the number $0.a_1a_1a_2a_2a_3a_3\dots$ will satisfy Theorem \ref{thm:one}. See Remark \ref{rem:trivial}.

This, in turn, leads to a deeper question: if a number being normal along non-trivial arithmetic progressions is not enough to guarantee normality of the original number, are there other non-trivial sequences one could select along which would, collectively, imply normality?

First let us consider which sequences trivially give normality. The following result 
(whose first half is well-known) says that a sequence trivially implies normality if and only if the sequence has asymptotic lower density equal to $1$.  The asymptotic lower density of an increasing sequence $A=\{n_1,n_2,n_3,\dots\}\subset \mathbb{N}$ is equal to $\liminf_{N\to \infty} |A\cap [1,N]|/N$.

\begin{prop}\label{prop:trivial}
Let $n_1<n_2<n_3<\dots$ be an increasing sequence of natural numbers. 

If $0.a_{n_1}a_{n_2}a_{n_3}\dots\in[0,1)$ is normal and  the asymptotic lower density of the sequence of $n_i$'s is equal to $1$, then $0.a_1a_2a_3\dots$ is normal.

On the other hand, if the asymptotic lower density of the sequence of $n_i$'s is strictly less than $1$, then there exist numbers $0.a_1a_2a_3\dots\in[0,1)$ which are not normal, even though $0.a_{n_1}a_{n_2}a_{n_3}\dots$ is normal.
\end{prop}

By altering the method of proving this proposition, we can show a condition by which normality along non-trivial sequences does imply normality overall. In particular if we have a collection of increasing sequences whose asymptotic lower density converges to $1$, then normality along these sequences implies normality overall.

\begin{thm}\label{thm:two}
Let $0.a_1a_2a_3\dots\in [0,1) $ and suppose that for any $\epsilon>0$ there exists an increasing sequence $n_1<n_2<n_3<\dots$ of positive integers with asymptotic lower density greater than $1-\epsilon$ such that $0.a_{n_1}a_{n_2}a_{n_3}\dots$ is normal.  Then $0.a_1a_2a_3\dots$ is normal.
\end{thm}

In the next result we will show that Theorem \ref{thm:two} is close to being sharp. For the purposes of this result, a set $A\subset \mathbb{N}$ is periodic if there exists a $m\in \mathbb{N}$ such that $(A-m)\cap \mathbb{N}=A$. Any such $m$ satisfying this condition will be called a period of $\mathbb{N}$.

\begin{thm} \label{thm:three}
Let $\mathcal{N}$ be a collection of periodic increasing sequences $n_1<n_2<n_3<\dots$ of the positive integers. Suppose there exists $K,L\in \mathbb{N}$ such that $\{K(n-1)+L,K(n-1)+L+1,K(n-1)+L+2,\dots,Kn+L-1\}$ is not a subset of any of the sequences for any $n\in \mathbb{N}$.

Then there exists a real number $0.a_1a_2a_3\dots\in [0,1)$ that is not normal, yet $0.a_{n_1}a_{n_2}a_{n_3}\dots$ is normal for all sequences in $\mathcal{N}$.
\end{thm}

In particular, the condition applied to $\mathcal{N}$ guarantees it cannot contain a periodic sequence of density greater than $1-1/L$. However, the additional restriction on sequences in $\mathcal{N}$ is a question of the thickness of a subset of $\mathbb{N}$, here referring to the length of allowable sub-sequences  of consecutive integers. So this leaves open a question of whether the condition on asymptotic lower density is the right one, or whether we should be using a condition on thickness instead.

We conclude the introduction by noting that we seem to have come very far afield from Wall's theorem. Wall's theorem states that selecting along an arithmetic progression presesrves normality, but selecting along an arbitrary increasing sequence may not preserve normality. It was shown by Kamae and Weiss \cite{Kamae,Weiss} that a sequence is guaranteed to preserve normality if and only if it is ``deterministic" and has positive asymptotic lower density. One could think of this as a generalized Wall's theorem.

We will come back to the proper definition of deterministic later, as it is complicated. Here we only note that determistic sequences are a subset of all sequences, so we could add the requirement in Theorem \ref{thm:two} that all sequences under consideration are deterministic. Periodic sequences like those in Theorem \ref{thm:three} are a type of deterministic sequence, and we could likely weaken the condition of periodicity to a condition of determinism, and in this sense could see that the altered Theorem \ref{thm:two} does come close to being a sharp converse to the generalized Wall's theorem. However, in the interest of keeping the paper short and readable, we will not give the proof of this.

\section{Preliminaries}
\subsection{Strings of strings}

Given a finite set of digits $\mathcal{D}$, we let the set of strings of length $k$ to be an ordered $k$-tuple with elements in $\mathcal{D}$, and denote this by $\mathcal{D}^k$ in the usual way. We will often write $[d_1,d_2,\dots,d_k]$ for such an ordered $k$-tuple.

We may then compose this notation and write, for example, $(\mathcal{D}^k)^\ell$. By this we mean the set of all ordered $\ell$-tuples whose elements belong to the set $\mathcal{D}^k$ of ordered $k$-tuples. Such an element might look like
\[
\left[ [d_1,d_2,\dots,d_k],[d_{k+1},d_{k+2},\dots,d_{2k}],\dots,[d_{k(\ell-1)+1},\dots,d_{k\ell}]\right].
\]
There is a standard bijection from such elements to $\mathcal{D}^{k\ell}$, and the element above would be mapped to
\[
[d_1,d_2,\dots,d_{k\ell}].
\]
For the rest of this paper, whenever we refer to considering or interpreting an element of $(\mathcal{D}^k)^\ell$ as an element of $\mathcal{D}^{k\ell}$ (or vice-versa), we mean that we are applying this standard bijection.

Note that we may allow $\ell=\infty$ and as such would have a standard bijection between $(\mathcal{D}^k)^\infty$ and $\mathcal{D}^\infty$. For slightly easier readibility, if we have an infinite tuple, we will use regular parenthesis $(\cdot )$ rather than brackets $[\cdot ]$.

\subsection{Symbolic Bernoulli shifts}

While we could express all our results merely in terms of real numbers, as the last section hints it will be easier for us if we instead treat them as manipulations of infinite strings of digits. In this section we will go over the basics of how to do this.

Let $\mathcal{D}$ be a finite set of digits and let $X=\mathcal{D}^\infty$. We will describe points $x\in X$ by \[ x=(a_1(x),a_2(x),a_3(x),\dots)=(a_1,a_2,a_3,\dots)\] with each $a_i\in \mathcal{D}$.

For a finite string $s=[d_1,d_2,\dots,d_k]\in\mathcal{D}^k$ we define the cylinder sets $C_s$ to be all elements $x\in X$ such that $a_1(x)=d_1, a_2(x)=d_2,\dots,a_k(x)=d_k$. 

For each $d\in \mathcal{D}$, let $\lambda_d$ be a non-negative number such that $\sum_{d\in \mathcal{D}} \lambda_d=1$. Then for each finite string $s=[d_1,d_2,\dots,d_k]$, we define $\mu(C_s)$ to be $\prod_{i=1}^k \lambda_{d_i}$. We use the cylinder sets to generate a $\sigma$-algebra and extend $\mu$ to be a measure on this $\sigma$-algebra.

Finally, we let $T$ be the standard forward shift on this space. So $T(a_1,a_2,a_3,\dots)=(a_2,a_3,a_4,\dots)$. We will refer to the dynamical system $(X,\mu,T)$ as a Bernoulli shift on the digit set $\mathcal{D}$.

We say that a point $x=(a_1,a_2,a_3,\dots)\in X$ is normal with respect to this transformation if for all finite strings $s=[d_1,d_2,\dots,d_k]$, we have
\[
\lim_{n\to \infty} \frac{\#\{0\le i \le n-1: a_{i+j}=d_j,j=1,2,\dots,k\}}{n} = \mu(C_s).
\]
This limit can be rephrased in a more standard ergodic fashion as 
\[
\lim_{n\to \infty} \frac{\#\{0\le i \le n-1:T^i x\in C_s\}}{n} = \mu(C_s).
\]

Consider the Bernoulli shift $(X,\mu,T)$ on the digit set $\mathcal{D}=\{0,1,\dots,b-1\}$ with $\mu$ defined by $\lambda_d=1/b$ for all $d\in \mathcal{D}$. Then this is clearly a symbolic representation of a base-$b$ expansion, with a natural correspondence given by $(a_1,a_2,a_3,\dots)\leftrightarrow 0.a_1a_2a_3\dots$. (This is well-defined up to a measure-zero set that can be ignored for the purposes of this paper.) This correspondence also clearly preserves  normality. We will therefore, for the rest of this paper, consider all base-$b$ systems as Bernoulli shifts.

In a more general setting, such symbolic shifts correspond to generalized L\"{u}roth series. As such, although normal points $x\in X$ have full measure by Birkhoff's pointwise ergodic theorem, if an explicit construction of such a point is desired, then examples can be found in \cite{AP,MM,VLuroth}.

\begin{rem}
All of the theorems given in the introduction are stated with respect to the base-$b$ expansion, they all hold for any Bernoulli shift. This is because there is no point where we make special use that the measure of $C_s$ for a length-$k$ string $s$ is $b^{-k}$. We only make use of the fact that the measure is a product measure on the digits.

In fact, we could allow $\mathcal{D}$ to be countably infinite and all the results would still hold.

However, for ease of readibility, we will express all the proofs with respect to the base-$b$ expansion given in the introduction.
\end{rem}

\subsection{Normality with respect to $T$ and $T^k$}

Let $(X,\mu,T)$ be a Bernoulli shift on the digit set $\mathcal{D}$ as defined above. Consider the Bernoulli shift $(X_k,\mu_k,T^k)$ with $X_k = (\mathcal{D}^k)^\infty$ and $\mu_k$ is defined via $\lambda_{[d_1,d_2,\dots,d_k]}=\mu(C_{[d_1,d_2,\dots,d_k]})$. Since there is a natural bijection between $X$ and $X_k$, we refer to the single forward shift on $X_k$ by $T^k$, since it is acting by $T^k$ on $X$.

It makes sense to refer to a point $x\in X$ as also belonging to $X_k$, since we may apply the standard bijection to achieve the corresponding point in $X_k$.

\begin{lem}\label{lem:T^k}
Under the definitions above, a point $x\in X$ is normal with respect to $(X,\mu,T)$ if and only if it is normal with respect to $(X_k,\mu_k,T^k)$ when seen as an element of $X_k$.
\end{lem}

Schweiger \cite{Schweiger} was the first to state this result, although his proof in one direction was erroneous. See \cite{VJoint} for a corrected proof. In the special case of base-$b$ expansions, this was found several years earlier. See \cite{Schmidt}.

\subsection{Deterministic sequences}\label{sec:deterministic}
In the introduction, we briefly made mention of deterministic sequences. For completeness, let us define better what we mean.

Consider an increasing sequence of positive integers $n_1<n_2<n_3<\dots$ and let $\omega\in \{0,1\}^\infty$ be such that $\omega_n=1$ if and only if $n=n_i$ for some $i$. A sequence is said to be completely deterministic in the sense of Weiss if all the weak-limits of the set of empirical measures for the forward shifts of $\omega$ have zero measure-theoretic entropy. Rauzy \cite{Rauzy} provided an alternative definition where a sequence is deterministic if 
\[
\lim_{s\to \infty} \limsup_{N\to \infty} \inf_{\phi\in E_s} \frac{1}{N} \sum_{n<N} \min\{1,|\omega_n-\phi(\omega_{n+1},\dots,\omega_{n+s})|\} = 0,
\]
where $E_s$ is the set of all functions from $\{0,1,\dots,b-1\}^s$ to $\{0,1,\dots,b-1\}$. One can think of the function $\phi$ as an attempt to guess at the value of $\omega_n$ given knowledge of $\omega_{n+1},\dots,\omega_{n+s}$. So Rauzy's definition says that a sequence is deterministic if the value of $\omega_n$ is ``determined" by the tail $\omega_{n+1},\omega_{n+2},\dots$.

That all periodic sequences are deterministic follows from either definition. However, we will make use of a special case of a result of Auslander and Dowker \cite[Theorem~6]{AD}, as it sets up the connection to normality most clearly.

\begin{prop}\label{prop:AD}
Let $(Y,\mathcal{G},\nu)$ be a compact measure space with $\nu(Y)=1$. Let $S:Y\to Y$ be a $\nu$-measure-preserving invertible transformation that has zero entropy.
 Let $U\subset Y$ be an open set with $\nu(U)>0$ and $\nu(\partial U)>0$. Let $y_0\in Y$ be generic.

Let $n_1<n_2<n_3<\dots$ be an increasing sequence of positive integers such that $n=n_i$ for some $i$ if and only if $T^n y_0\in U$. In other words, the $n_i$'s are the sequence of visiting times for the orbit of $y_0$ to the set $U$.

Then if $(a_1,a_2,a_3,\dots)$ is normal with respect to some Bernoulli shift, then $(a_{n_1},a_{n_2},a_{n_3}, \\ \dots)$ is also normal with respect to the same Bernoulli shift.
\end{prop}

We remark that Auslander and Dowker technically only proved this the standard base-$2$ Bernoulli shift, but it is a simple tweak of their proof to get the result above.

Suppose for some positive integer $m\ge 2$, $Y=\{0,1,\dots,m-1\}$, $\nu$ is the normalized counting measure on $X$, $S$ is given by $Sx=x+1\pmod{m}$, $U$ is any subset of $Y$, and $y_0$ is any element of $Y$. This can be used so that the corresponding $n_i$'s are any desired periodic sequence and so gives the following as an immediate consequence.

\begin{lem}
Let $(a_1,a_2,a_3,\dots)$ be normal with respect to some Bernoulli shift $(X,\mu,T)$ and let $n_1<n_2<n_3<\dots$ be any eventually periodic sequence. Then $(a_{n_1},a_{n_2},a_{n_3}, \dots)$ is also normal with respect to $(X,\mu,T)$.
\end{lem}

The sequences covered by Auslander and Dowker's result are quite varied. For instance, they cover generalized linear functions such as $n_i=[\alpha i +\beta]$ for $\alpha>1,\beta\ge 0$. However, it is not clear whether they cover all possible deterministic sequences.

\section{Proof of Proposition \ref{prop:trivial}}

Let $x=(a_1,a_2,a_3,\dots)$ belong to the Bernoulli base-$b$ shift.

Let $n_1<n_2<n_3<\dots$ be an increasing sequence of natural numbers with asymptotic lower density equal to $1$. Suppose that $y=(a_{n_1},a_{n_2},a_{n_3},\dots)$ is normal with respect to this same Bernoulli shift. 

Consider an arbitrary string $s=[d_1,d_2,\dots,d_k]$ with digits belonging to the digit set $\{0,1,\dots,b-1\}$. We wish to show that the limiting frequency of $s$ in $x$ is $b^{-k}$. 

Let $N$ be a large positive integer. Let $j=j(N)$ denote the largest index such that $n_j\le N$. Since the asymptotic lower density of the sequence is $1$, we have that $j(N)=N(1+o(1))$.

So consider the number of times that $s$ appears starting in the first $N$ digits of $x$. Each such string will also appear in the first $j(N)$ digits of $y$ unless one of the digits of the string gets removed in going from $x$ to $y$. This happens at most $o(N)$ times. Similarly any such string appearing in the first $j(N)$ digits of $y$ appears in the first $N$ digits of $x$ unless a digit was inserted somewhere in the middle of it, which happens again at most $o(N)$ times.

Thus, we have that
\begin{align*}
\frac{\#\{0\le i \le N-1: T^i x\in C_s\}}{N} &= \frac{\#\{0\le i \le j(N)-1:T^i y \in C_s\}+o(N)}{N}\\
&=  \frac{\#\{0\le i \le j(N)-1:T^i y \in C_s\}}{j(N)}(1+o(1))\\
&= b^{-k}(1+o(1)),
\end{align*}
by the normality of $y$. Thus $x$ is normal.

For the second part of the proposition, suppose $n_1<n_2<n_3<\dots$ is an increasing sequence of natural numbers with asymptotic lower density $\alpha <1$. Suppose that $y=(a_{n_1},a_{n_2},a_{n_3},\dots)$ is normal.

Let $x=(a_1,a_2,a_3,\dots)$ be defined so that $a_n=0$ if $n\neq n_i$ for any $i$. 

Let $N$ be an integer such that $j(N)$ (as defined above) is at most $N(1+\alpha)/2$. By our assumption of the density of the sequence, there must be arbitrarily large such $N$'s.

For such an $N$, consider how many $0$'s appear in the first $N$ digits of $x$. By the normality of $y$, there must be $j(N)b^{-k}(1+o(1))$ such $0$'s coming from the $0$'s of $y$, and there are also $N-j(N)$ such $0$'s coming from the digits $a_n$ with $n=n_i$ for any $i$.

Thus the total number of $0$'s in the first $N$ digits of $x$ is
\begin{align*}
(N-j(N))+j(N)b^{-k}(1+o(1))&=Nb^{-k}+(N-j(N))+(j(N)-N)b^{-k} +o(j(N))\\
&= Nb^{-k}+(N-j(N))(1-b^{-k})+o(N)\\ &\ge N\left( b^{-k} +\left(1-\frac{1+\alpha}{2}\right)\left(1-b^{-k}\right)+o(1)\right).
\end{align*}
And since, for certain arbitrarily large $N$, this will exceed $N(b^{-k}+\epsilon)$ for some small $\epsilon>0$, we have that $x$ is not normal.

\section{Proof of Theorem \ref{thm:two}}

Let $x=(a_1,a_2,a_3,\dots)$ satisfy the conditions of the theorem. 

Consider an arbitrary finite string $s=[d_1,d_2,\dots,d_k]$. We want to show that the limiting frequency of $s$ in $x$ is $b^{-k}$.

Select $\epsilon>0$ arbitrary and pick an increasing sequence of positive integers, $n_1<n_2<n_3<\dots$, whose asymptotic lower density strictly exceeds $1-\epsilon$, such that $y=(a_{n_1},a_{n_2},a_{n_3},\dots)$ is normal.

Now we borrow several ideas from the proof of Proposition \ref{prop:trivial}. First let $N$ and $j(N)$ be defined as in that proof. We will assume that $N$ is sufficiently large so that $j(N)\ge N(1-\epsilon)$; in particular $j(n)=N(1+O(\epsilon))$. Then by the argument of the previous proof, we have that the number of times $s$ appears starting in the first $N$ digits of $x$ is equal to the number of times $s$ appears in the first $j(N)$ digits of $y$, up to an error of $O(k\epsilon N)$. (The presence of $k$ in this term comes from the fact that, for example, if we delete a single digit from $x$, this alters $k$ different strings of length $k$.)

Thus, again mimicking the previous proof, we have
\begin{align*}
\frac{\#\{0\le i \le N-1:T^i x\in C_s\}}{N} &= \frac{\#\{0\le i \le j(N)-1: T^i y\in C_s\} +O(k\epsilon N)}{j(N)(1+O(\epsilon))}\\
&=  \frac{\#\{0\le i \le j(N)-1: T^i y\in C_s\} }{j(N)}\cdot \frac{1+O(k\epsilon)}{1+O(\epsilon)}\\
&= b^{-k}\cdot \frac{(1+O(k\epsilon))(1+o(1))}{1+O(\epsilon)}.
\end{align*}
Now by letting $N$ go to infinty, we get that the limiting frequency of $s$ is $b^{-k}(1+O(k\epsilon))/(1+O(\epsilon))$. Then, by taking $\epsilon$ arbitrarily small, we get the desired limiting frequency of $b^{-k}$.

\section{Proof of Theorem \ref{thm:three}}

Let $\mathcal{N},K,L$ be as in the statement of the theorem. Without loss of generality, we may assume $L=1$. Let $(X,\mu,T)$ be the usual base-$b$ symbolic shift.

Let us consider a new symbolic shift on $\mathcal{D}^K$, where $\mathcal{D}=\{0,1,2,\dots,b-1\}$. We will define the measure $\nu$ on $Y=(\mathcal{D}^K)^\infty$ by
\begin{equation}\label{eq:badcasedef}
\lambda_{[d_1,d_2,\dots,d_K]}= \begin{cases}
\dfrac{1}{b^K} -\dfrac{(-1)^{d_1+\dots+d_K}}{2b^K} & \text{ if all the }d_i\text{ are either } 0 \text{ or }1\\
\dfrac{1}{b^K} & \text{ otherwise.}
\end{cases}
\end{equation}
We let the forward shift on this space be called $T_Y$

Now we wish to consider what we will call ``starred digits," $\mathcal{S}=(\mathcal{D}\cup\{*\})^K \setminus \mathcal{D}^K$. These can be seen as elements of $\mathcal{D}^K$ with at least one digit replaced with $*$. We can then consider starred strings to be elements of $\mathcal{S}^m$ for some integer $m\ge 1$.

Although starred strings are defined over a larger digit set that includes $*$, we may interpret them as a collection of strings with digits in $\mathcal{D}^K$. In particular, a starred string in $\mathcal{S}^m$ can be considered as the collection of all strings in $(\mathcal{D}^K)^m$ where each $*$ in any of the starred digits is allowed to be replaced by any of the digits in $\mathcal{D}$. And, it should be emphasized, we don't have to use the same digit from $\mathcal{D}$ each time we do this replacement.

Thus, with this new interpretation, we may talk about a starred string $s\in\mathcal{S}^m$ ``appearing" in the expansion of a point in $(\mathcal{D}^K)^\infty$. In particular, $s$ appears in this point if one of the corresponding strings in $(\mathcal{D}^K)^m$ appears in this point. Likewise we may define the measure $\nu(C_s)$ to be the union of the $\nu$-measure of all the cylinder sets for the corresponding strings in $\mathcal{D}^K$. We then note that the relative frequency with which $s$ appears in a normal point equals the measure $\nu(C_s)$.

We claim that for any string $s\in \mathcal{S}^m$, we have $\nu(C_s) = b^{-n}$ where $n$ is the number of digits from $\mathcal{D}$ that appear in $s$ (when viewed as a string with $mK$ total digits from the set $\mathcal{D}\cup\{*\}$).

As an easy first case, consider $s=[[d_1,\dots,d_{K-1},*]]\in\mathcal{S}^1$. 
In this case we have
\[
\nu(C_s)=\sum_{d_K\in\mathcal{D}} \lambda_{[d_1,\dots,d_K]} = \frac{1}{b^{K-1}}.
\]
This follows because if $d_1,\dots,d_{K-1}$ are all either $0$ or $1$, then all of the summands are $b^{-K}$ except for one term of $1.5*b^{-K}$ and one term of $.5*b^{-K}$, and if at least one of the $d_1,\dots,d_{K-1}$ is not $0$ or $1$, then all of the summands are $b^{-K}$.

It is clear the same will hold for any $s\in\mathcal{S}^1$ that has only one $*$ in it.

Now let us consider, for example, $s=[[d_1,\dots,d_{K-2},*,*]]\in\mathcal{S}^1$.
Then we have
\[
\nu(C_s)=\sum_{d_{K-1},d_K\in\mathcal{D}} \lambda_{[d_1,\dots,d_K]} = \sum_{d_{K-1}\in \mathcal{D}} \frac{1}{b^{K-1}} = \frac{1}{b^{K-2}}.
\]
A similar result can be seen to hold for any $s\in\mathcal{S}^1$ that has exactly two $*$'s in it

The case of two $*$'s is very instructive, and from it we can clearly see that by induction $\nu(C_s)=b^{-n}$ for all $s\in \mathcal{S}^1$. And since $(Y,\nu,T_Y)$ as defined in this section is Bernoulli, we see that $\nu(C_s)=b^{-n}$ for all $s\in\mathcal{S}^m$ for any integer $m\ge 1$.

\begin{rem}\label{rem:trivial}
As an aside, note that if $K=2$ and we let $\lambda_{[d_1,d_2]}=1/b$ if $d_1=d_2$ and $0$ otherwise, then the measure of the starred strings is again given by $\nu(C_s)=b^{-n}$. This is why the construction given after Theorem \ref{thm:one} works.
\end{rem}

Let $x=([a_1,\dots,a_K],[a_{K+1},\dots,a_{2K}],\dots)$ be a normal point for this symbolic shift $(Y,\nu,T_Y)$, and let us apply the standard bijection to interpret it as a point in $\mathcal{D}^\infty$.

First we claim that $x$ cannot be normal with respect to $(X,\mu,T)$. If it were, then by Lemma \ref{lem:T^k} it would be normal with respect to $(X_K,\mu_K,T^K)$. In particular the limiting frequency of the digit $[0,0,\dots,0]\in \mathcal{D}^K$ should be $b^{-K}$. However, by construction, the actual limiting frequency of the digit $[0,0,\dots,0]\in\mathcal{D}^K$ is $\nu(C_{[[0,\dots,0]]})$, which will either be $.5b^{-K}$ or $1.5b^{-K}$. So $x$ cannot be normal with respect to $(X,\mu,T)$.

On the other hand, we claim that after selecting along any of the periodic sequences in $\mathcal{N}$, $x$ is normal. In particular, suppose that we are looking at an increasing periodic sequence of positive integers $n_1<n_2<n_3<\dots$ with period $m$. Then this sequence also has period $mK$.  Suppose the periodic sequence contains $p$ elements in the interval $[1,mK]$---in other words, assume that $n_1<n_2<\dots<n_p\le mK < n_{p+1}$.

Let $y=(a_{n_1},a_{n_2},a_{n_3},\dots)\in \mathcal{D}^\infty$. Then this is normal, by lemma \ref{lem:T^k}, if and only if $y_p=([a_{n_1},\dots,a_{n_p}],[a_{n_{p+1}},\dots,a_{n_{2p}}],\dots)$ is normal in $(X_p,\mu_p,T^p)$. Consider any string $s_p=[\mathfrak{d}_1,\mathfrak{d}_2,\dots,\mathfrak{d}_{pj}]\in (\mathcal{D}^p)^j$ of length $j\ge1$. To complete the proof, we wish to show that the limiting frequency that $s_p$ occurs in $y_p$ equals $b^{-pj}$.

Now consider a  string $s=[[d_1,\dots,d_K],\dots,[d_{(mj-1)K+1},\dots,d_{mjK}]]\in ((\mathcal{D}\cup\{*\})^K)^{mj}$ defined in the following way:
\[
d_i =\begin{cases}
\mathfrak{d}_{i'}, & \text{ if }i=n_{i'} \text{ for some }i'\in \mathbb{N},\\
*, & \text{ otherwise}.
\end{cases}
\]
The restriction placed on the sequences in $\mathcal{N}$ guarantees that for each set $\{K(n-1)+1,K(n-1)+2,\dots,Kn\}$, at least one of these elements does not belong to the sequence of $n_i$'s and so each of the $mj$ digits of $s$ contains at least one $*$. Thus, $s$ truly is an element of $\mathcal{S}^{mj}$, not just $((\mathcal{D}\cup\{*\})^K)^{mj}$. 

Now by the work we did earlier, the limiting frequency with which $s_p$ occurs in $y_p$ is equal to the frequency with which the starred string $s$ occurs in $x_K$; however, this is exactly $\nu(C_s)=b^{-pj}$ as desired.

This completes the proof.

\section{Acknowledgments}

The author would like to thank Bill Mance for inspiring this line of questioning and Vitaly Bergelson for offering helpful questions.


\begin{thebibliography}{10}

\bibitem{AP}
Max Aehle and Matthias Paulsen.
\newblock Construction of normal numbers with respect to generalized l\"{u}roth
  series from equidistributed sequences.
\newblock {\em arXiv preprint arXiv:1509.08345}, 2015.

\bibitem{AD}
J.~Auslander and Y.~N. Dowker.
\newblock On disjointness of dynamical systems.
\newblock {\em Math. Proc. Cambridge Philos. Soc.}, 85(3):477--491, 1979.

\bibitem{Kamae}
Teturo Kamae.
\newblock Subsequences of normal sequences.
\newblock {\em Israel J. Math.}, 16:121--149, 1973.

\bibitem{MM}
Manfred~G. Madritsch and Bill Mance.
\newblock Construction of {$\mu$}-normal sequences.
\newblock {\em Monatsh. Math.}, 179(2):259--280, 2016.

\bibitem{Rauzy}
G.~Rauzy.
\newblock Nombres normaux et processus d\'eterministes.
\newblock {\em Acta Arith.}, 29(3):211--225, 1976.

\bibitem{Schmidt}
Wolfgang~M. Schmidt.
\newblock On normal numbers.
\newblock {\em Pacific J. Math.}, 10:661--672, 1960.

\bibitem{Schweiger}
F.~Schweiger.
\newblock Normalit\"at bez\"uglich zahlentheoretischer {T}ransformationen.
\newblock {\em J. Number Theory}, 1:390--397, 1969.

\bibitem{VLuroth}
J.~Vandehey.
\newblock A simpler normal number construction for simple {L}\"uroth series.
\newblock {\em J. Integer Seq.}, 17(6):Article 14.6.1, 18, 2014.

\bibitem{VJoint}
Joseph Vandehey.
\newblock On the joint normality of certain digit expansions.
\newblock {\em arXiv preprint arXiv:1408.0435}, 2014.

\bibitem{WallThesis}
Donald~D. Wall.
\newblock {\em N{ORMAL} {NUMBERS}}.
\newblock ProQuest LLC, Ann Arbor, MI, 1950.
\newblock Thesis (Ph.D.)--University of California, Berkeley.

\bibitem{Weiss}
B~Weiss.
\newblock Normal sequences as collectives.
\newblock In {\em Proc. Symp. on Topological Dynamics and Ergodic Theory, Univ.
  of Kentucky}, 1971.

\end{thebibliography}
\end{document}